\documentclass{article}
\usepackage{amssymb,amsmath}
\newtheorem{thm}{Thm.}[section]
\newtheorem{def.}{Def.}[section]
\newtheorem{rem}{Remark}[section]
\newtheorem{prop}{Prop.}[section]
\newtheorem{cor}{Cor.}[section]

\input xy.sty
\xyoption{all}
\xyoption{graph}
\xyoption{poly}
\xyoption{knot}
\xyoption{rotate}
\begin{document}

\def\TC{\xygraph{!{0;/r.75pc/:}
!P8"c"{~:{(3.5,0):}~>{}}
!{\ar @{->} "c3";"c6",^r}}}
\def\TD{\xygraph{!{0;/r.75pc/:}
!P4"a"{~:{(1.75,0):}~>{}}!P8"b"{~:{(3.5,0):}~>{}}
!{\hunder~{"a2"}{"a1"}{"a3"}{"a4"}}
!{\ar @{<-} "a2";"b3",^{r'_1}}
!{\ar @{->} "a3";"b6",^{r'_2}}
!{\save0;"a3"-"a4":"a4",
\hcap-@(+.1)\restore}}}
\def\TA1{\xygraph{!{0;/r.75pc/:}
!P8"c"{~:{(3.5,0):}~>{}}
!{\ar @{->} "c3";"c6",^r}}}
\def\TB1{\xygraph{!{0;/r.75pc/:}
!P4"a"{~:{(1.75,0):}~>{}}!P8"b"{~:{(3.5,0):}~>{}}
!{\hover~{"a2"}{"a1"}{"a3"}{"a4"}}
!{\ar @{<-} "a2";"b3",^{r'_1}}
!{\ar @{->} "a3";"b6",^{r'_2}}
!{\save0;"a3"-"a4":"a4",
\hcap-@(+.1)\restore}}}
\def\RA{\xygraph{!{0;/r.6pc/:}
!P4"a"{~:{(-.75,-.75):}~>{}}!P8"b"{~:{(6,0):}~>{}}
!{\POS"b3" \ar @/^3ex/ "b6", _{r}}
!{\POS"b2" \ar @/_3ex/ "b7", ^{s}}}}
\def\RB{\xygraph{!{0;/r.6pc/:}
!P4"a"{~:{(-.5,-.5):}~>{}}!P8"b"{~:{(6,0):}~>{}}
!{\xunderv~{"b2"}{"b3"}{"a2"}{"a4"}@(.65)<>|><{r'_1}|{s'}}
!{\xunderv~{"a4"}{"a2"}{"b6"}{"b7"}@(.35)<{r'_2}>{r'_3}}}}
\def\QA{\xygraph{!{0;/r.6pc/:}
!P4"a"{~:{(-.75,-.75):}~>{}}!P8"b"{~:{(6,0):}~>{}}
!{\POS"b3" \ar @/^3ex/ "b6", _{r}}
!{\POS"b7" \ar @/^3ex/ "b2", _{s}}}}
\def\QB{\xygraph{!{0;/r.6pc/:}
!P4"a"{~:{(-.5,-.5):}~>{}}!P8"b"{~:{(6,0):}~>{}}
!{\xunderv~{"b2"}{"b3"}{"a2"}{"a4"}@(.65)<>|<<{r'_1}|{s'}}
!{\xunderv~{"a4"}{"a2"}{"b6"}{"b7"}@(.35)<{r'_2}>{r'_3}}}}
\def\Reid{
\xy 0;/r1pc/:
\xoverv~{(-1,3)}{(1,3)}{(-1,1)}{(1,1)}|<<{r'}|{s'}>{u'}  \xoverv~{(-4,-2)}{(-4,-.5)}{(0,-.5)}{(-1,1)}<{v'}|{t'}
\xunderv~{(4,-.5)}{(4,-2)}{(1,1)}{(0,-.5)}|<<{w'}
\endxy}

\def\dieR{
\xy 0;/r-1pc/:
\xoverv~{(-1,3)}{(1,3)}{(-1,1)}{(1,1)}<{w}|{v}>{u}  \xoverv~{(-4,-2)}{(-4,-.5)}{(0,-.5)}{(-1,1)}<<|<<{s}|{t}
\xunderv~{(4,-.5)}{(4,-2)}{(1,1)}{(0,-.5)}<{r}
\endxy}

\newread\epsffilein    
\newif\ifepsffileok    
\newif\ifepsfbbfound   
\newif\ifepsfverbose   
\newif\ifepsfdraft     
\newdimen\epsfxsize    
\newdimen\epsfysize    
\newdimen\epsftsize    
\newdimen\epsfrsize    
\newdimen\epsftmp      
\newdimen\pspoints     
\pspoints=1bp          
\epsfxsize=0pt         
\epsfysize=0pt         
\def\epsfbox#1{\global\def\epsfllx{72}\global\def\epsflly{72}%
   \global\def\epsfurx{540}\global\def\epsfury{720}%
   \def\lbracket{[}\def\testit{#1}\ifx\testit\lbracket
   \let\next=\epsfgetlitbb\else\let\next=\epsfnormal\fi\next{#1}}%
\def\epsfgetlitbb#1#2 #3 #4 #5]#6{\epsfgrab #2 #3 #4 #5 .\\%
   \epsfsetgraph{#6}}%
\def\epsfnormal#1{\epsfgetbb{#1}\epsfsetgraph{#1}}%
\def\epsfgetbb#1{%
%
%
\openin\epsffilein=#1
\ifeof\epsffilein\errmessage{I couldn't open #1, will ignore it}\else
%
%
   {\epsffileoktrue \chardef\other=12
    \def\do##1{\catcode`##1=\other}\dospecials \catcode`\ =10
    \loop
       \read\epsffilein to \epsffileline
       \ifeof\epsffilein\epsffileokfalse\else
%
%
          \expandafter\epsfaux\epsffileline:. \\%
       \fi
   \ifepsffileok\repeat
   \ifepsfbbfound\else
    \ifepsfverbose\message{No bounding box comment in #1; using defaults}\fi\fi
   }\closein\epsffilein\fi}%
%
%
\def\epsfclipon{\def\epsfclipstring{ clip}}%
\def\epsfclipoff{\def\epsfclipstring{\ifepsfdraft\space clip\fi}}%
\epsfclipoff
\def\epsfsetgraph#1{%
   \epsfrsize=\epsfury\pspoints
   \advance\epsfrsize by-\epsflly\pspoints
   \epsftsize=\epsfurx\pspoints
   \advance\epsftsize by-\epsfllx\pspoints
%
%
   \epsfxsize\epsfsize\epsftsize\epsfrsize
   \ifnum\epsfxsize=0 \ifnum\epsfysize=0
      \epsfxsize=\epsftsize \epsfysize=\epsfrsize
      \epsfrsize=0pt
%
%
     \else\epsftmp=\epsftsize \divide\epsftmp\epsfrsize
       \epsfxsize=\epsfysize \multiply\epsfxsize\epsftmp
       \multiply\epsftmp\epsfrsize \advance\epsftsize-\epsftmp
       \epsftmp=\epsfysize
       \loop \advance\epsftsize\epsftsize \divide\epsftmp 2
       \ifnum\epsftmp>0
          \ifnum\epsftsize<\epsfrsize\else
             \advance\epsftsize-\epsfrsize \advance\epsfxsize\epsftmp \fi
       \repeat
       \epsfrsize=0pt
     \fi
   \else \ifnum\epsfysize=0
     \epsftmp=\epsfrsize \divide\epsftmp\epsftsize
     \epsfysize=\epsfxsize \multiply\epsfysize\epsftmp   
     \multiply\epsftmp\epsftsize \advance\epsfrsize-\epsftmp
     \epsftmp=\epsfxsize
     \loop \advance\epsfrsize\epsfrsize \divide\epsftmp 2
     \ifnum\epsftmp>0
        \ifnum\epsfrsize<\epsftsize\else
           \advance\epsfrsize-\epsftsize \advance\epsfysize\epsftmp \fi
     \repeat
     \epsfrsize=0pt
    \else
     \epsfrsize=\epsfysize
    \fi
   \fi
%
%
   \ifepsfverbose\message{#1: width=\the\epsfxsize, height=\the\epsfysize}\fi
   \epsftmp=10\epsfxsize \divide\epsftmp\pspoints
   \vbox to\epsfysize{\vfil\hbox to\epsfxsize{%
      \ifnum\epsfrsize=0\relax
        \includegraphics{\ifepsfdraft}%
      \else
        \epsfrsize=10\epsfysize \divide\epsfrsize\pspoints
        \includegraphics{\ifepsfdraft}%
      \fi
      \hfil}}%
\global\epsfxsize=0pt\global\epsfysize=0pt}%
%
%
{\catcode`\%=12 \global\let\epsfpercent=
%
%
\long\def\epsfaux#1#2:#3\\{\ifx#1\epsfpercent
   \def\testit{#2}\ifx\testit\epsfbblit
      \epsfgrab #3 . . . \\%
      \epsffileokfalse
      \global\epsfbbfoundtrue
   \fi\else\ifx#1\par\else\epsffileokfalse\fi\fi}%
%
%
\def\epsfempty{}%
\def\epsfgrab #1 #2 #3 #4 #5\\{%
\global\def\epsfllx{#1}\ifx\epsfllx\epsfempty
      \epsfgrab #2 #3 #4 #5 .\\\else
   \global\def\epsflly{#2}%
   \global\def\epsfurx{#3}\global\def\epsfury{#4}\fi}%
%
%
\def\epsfsize#1#2{\epsfxsize}
%
%
\let\epsffile=\epsfbox

\title{{Quandles at Finite Temperatures I}}
		\author{Pedro Lopes\\
	        Departamento de Matem\'atica and Centro de Matem\'atica Aplicada\\
	        Instituto Superior Tecnico\\
	        Av. Rovisco Pais\\
		1049-001 Lisboa\\
	        Portugal\\
	        \texttt{pelopes@math.ist.utl.pt}}

\date{May 11, 2001}

\maketitle

\begin{abstract}
In \cite{jsCetal} quandle cohomology is used to produce invariants for particular embeddings of codimension two; 2-cocycles give rise to invariants for (classical) knots and 3-cocycles give rise to invariants for knotted surfaces. This is done by way of a notion of coloring of a diagram. Also, these invariants have the form of the state-sums (or partition functions) used in Statistical Mechanics.

By a careful analysis of these colorings of diagrams we are able to come up with new invariants which correspond to calculation of the partition function at finite temperatures.

\end{abstract}

\section{Introduction} \label{S:intro}
The first appearance of quandles in connection to knot theory seems to have been in the PhD thesis of D. Joyce - see his paper \cite{dJoyce} and also \cite{sMatveev}. One of the interesting notions there is understanding the (knot) quandle as a minimal algebraic structure which is invariant under Reidemeister moves. See also \cite{lhKauffman} for a discussion of further aspects of quandles.

In \cite{jsCetal} quandle cohomology (see \cite{rFenn1}, \cite{rFenn2}) is used to produce invariants for (classical) knots (see e.g. \cite{rhCrowell}), using 2-cocycles, and for knotted surfaces (see e.g. \cite{CetS}), using 3-cocycles. Important ingredients in the construction in \cite{jsCetal} are: the notion of a coloring of a diagram and the evaluation, for each coloring,  of products of cocycles (raised to $\pm 1$), which we shall here call ``Boltzmann weights''. It is then proved that the sum over all colorings of the Boltzmann weights yields an invariant. Thus, the invariants have the form of the state-sums (or partition functions) used in Statistical Mechanics with colorings standing for configurations (cf. \cite{hbCallen}).

Refining the ideas in \cite{jsCetal} we introduce a relation on the set of colorings which turns out to be an equivalence relation. Then we prove that the set of equivalence classes of colorings of a diagram is an invariant. Moreover, we prove that the set of Boltzmann weights associated to the set of equivalence classes of colorings of a diagram is also an invariant (by proving that equivalent colorings have the same Boltzmann weight). Hence, any symmetric function of the Boltzmann weights is an invariant. As a specialization of this last statement we recover the result in \cite{jsCetal} namely, that the sum of all Boltzmann weights is an invariant. Also, if we raise each Boltzmann weight to the same power and then add them we have again an invariant. So if one takes that power to be $-\frac{1}{kT}$ then the state-sum has the form of a partition function for a canonical ensemble at finite temperature.

\subsection{Organization} \label{S:org}

In section 2 the background material is introduced, i.e., quandles and quandle cohomology. In section 3 we refine the invariants for classical knots and in section 4 we refine the invariants for embedded surfaces. In section 5 we present ideas for future work.

\subsection{Acknowledgments} \label{S:ack}

This work was supported by the programme {\em Programa Operacional
``Ci\^{e}ncia, Tecnologia, Inova\c{c}\~{a}o''} (POCTI) of the
{\em Funda\c{c}\~{a}o para a Ci\^{e}ncia e a Tecnologia} (FCT),
cofinanced by the European Community fund FEDER.

It was done in IST, in Lisbon, Portugal, under supervision of Prof. L. Crane and Prof. R. F. Picken. The author would like to thank Prof. L. H. Kauffman for having introduced him to quandles and Prof. D. Roseman for the introduction to embeddings of surfaces in four space. Figures 1, 2, 3 and 4 are reproduced from \cite{jsCetal} by kind permission of the authors. The author is also grateful to Prof. J. S. Carter for the decisive boost for finishing Section 4 and for clearing a number of details. Thanks are also due to Prof. M. Viana and the staff at IMPA where the writing of this paper actually begun. A final special "thank you" to Paula for always being there...

\section{Background Material} \label{S:back}

In this section we introduce the material that we will make use of in the subsequent sections. Most of this is found in \cite{jsCetal} although we introduce further material from \cite{dJoyce}.

\begin{def.}[Quandle]A quandle is a set $X$ endowed with a binary operation, denoted $\triangleright$, such that:

(i) for any $a \in X$, $a\triangleright a = a$;

(ii) for any $a$ and $b \in X$, there is a unique $x \in X$ such that $a = x\triangleright b$;

(iii) for any $a$, $b$, and $c \in X$, $(a\triangleright b)\triangleright c = (a\triangleright c)\triangleright (b\triangleright c)$ (self-distributivity).
\end{def.} 

\begin{def.}[Rack]A rack is a set $X$ endowed with a binary operation that satisfies (ii) and (iii) above.\end{def.}

\vskip 5pt

The second axiom for quandles, above, implies the existence of a second binary operation in $X$, denoted $\triangleright ^{-1}$, defined as follows.

\begin{def.}[$\triangleright ^{-1}$]Let $X$ be a quandle. For any $a$, $b\in X$,

\[
a\triangleright ^{-1}b = x
\]
where $x$ is the unique element in $X$ such that $a = x\triangleright b$.\end{def.}

\begin{prop}If $\big( X,\triangleright \big)$ is a quandle then the following hold:

(i) for any $a\in X$, $a\triangleright ^{-1}a = a$;

(ii) for any $a$ and $b \in X$, $(a\triangleright ^{-1}b)\triangleright b = a = (a\triangleright b)\triangleright ^{-1}b$;

(iii) for any $a$, $b$, and $c \in X$, 
\begin{align}\notag
 (a\triangleright ^{-1}b)\triangleright ^{-1}c &= (a\triangleright ^{-1}c)\triangleright ^{-1}(b\triangleright ^{-1}c);\\ \notag
 (a\triangleright ^{-1}b)\triangleright c &= (a\triangleright c)\triangleright ^{-1}(b\triangleright c);\\ \notag
 (a\triangleright b)\triangleright ^{-1}c &= (a\triangleright ^{-1}c)\triangleright (b\triangleright ^{-1}c). 
\end{align}
\end{prop}Proof: (i) If $a\in X$,then $a\triangleright ^{-1}a$ is the unique element $x\in X$ such that $a=x\triangleright a$ which is $a$; 

(ii) for any $a$ and $b \in X$, $a\triangleright ^{-1}b$ is the unique element $x\in X$ such that $a=x\triangleright b$, hence $(a\triangleright ^{-1}b)\triangleright b = x\triangleright b = a$ and analogously for the other identity;

(iii) We will just prove the third equality; the others are dealt with in an analogous way. For any $a$, $b$, $c\in X$:

\[
(a\triangleright b)\triangleright ^{-1}c=w \Leftrightarrow w\triangleright c=a\triangleright b)
\]
Equivalently, 
\[
[(a\triangleright ^{-1}c)\triangleright (a\triangleright ^{-1}c)]\triangleright c=[(a\triangleright ^{-1}c)\triangleright c]\triangleright [(a\triangleright ^{-1}c)\triangleright c]=a\triangleright b
\]
and by unicity the equality follows.$\hfill \blacksquare$

\vskip 20pt

The basic example of a quandle is a group with $a\triangleright b := b^{-n}ab^{n}$, where $n$ is a fixed positive integer and juxtaposition denotes the group multiplication. Also, the axioms for quandles correspond to the three Reidemeister moves. In fact, an oriented knot diagram can be regarded as a table defining a quandle which is presented in the following way: the arcs stand for generators of the quandle and relations are to be read at each crossing as follows:

$$\xy 0;/r1pc/:
,{\xoverv[4.0] |> <{x}|{y}>{x\triangleright y}}
\endxy$$
i.e., if the over-arc (with orientation as indicated) is $y$ and the upper under-arc is $x$, then the lower under-arc is $x\triangleright y$ and, equivalently,

$$\xy 0;/r1pc/:
,{\xoverv[4.0] |< <{x}|{y}>{x{\triangleright ^{-1}y}}}
\endxy$$ 
i.e., if the over-arc (with orientation as indicated) is $y$ and the upper under-arc is $x$, then the lower under-arc is $x{\triangleright ^{-1}y}$.

In this way, Reidemeister equivalent diagrams give rise to different presentations of the same quandle, i.e., quandles so determined are knot invariants (cf. \cite{dJoyce}). This quandle is called the fundamental quandle of the knot.

\vskip 5pt

\begin{def.}For any quandle $X$, let $C_{n}^{R}(X)$ be the free abelian group generated by $n$-tuples $(x_1,...,x_n)$ of elements of the quandle $X$. Define a homomorphism $\partial _n : C_{n}^{R}(X) \rightarrow C_{n-1}^{R}(X)$

\begin{multline}
\partial (x_1,x_2,...,x_n) \\
= \sum _{i=2}^{n}(-1)^i[(x_1,x_2,...,x_{i-1},x_{i+1},...,x_n)\\
-(x_1*x_i,x_2*x_i,...,x_{i-1}*x_i,x_{i+1},...,x_n)]\notag
\end{multline}

for $n\ge 2$ and $\partial _n = 0$ for $n\le 1$.\end{def.}

\begin{prop}$C_{*}^{R}(X)=\{C_{n}^{R}(X),\partial _n\}$ is a chain complex.\end{prop}

Let $C_{n}^{D}(X)$ be the subset of $C_{n}^{R}(X)$ generated by $n$-tuples $(x_1,...,x_n)$ with $x_i=x_{i+1}$ for some $i\in \{1,2,...,n-1\}$ if $n\ge 2$; otherwise let $C_{n}^{D}(X)=0$. If $X$ is a quandle, then $\partial _n(C_{n}^{D}(X))\subset C_{n-1}^{D}(X)$ and $C_{*}^{D}(X)=\{C_{n}^{D}(X),\partial _n\}$ is a sub-complex of $C_{*}^{R}(X)$. Put $C_{n}^{Q}(X)=C_{n}^{R}(X)/C_{n}^{D}(X)$ and $C_{*}^{Q}(X)=\{C_{n}^{Q}(X),\partial' _n\}$ where $\partial' _n$ is the induced homomorphism. Henceforth, all boundary maps will be denoted by $\partial _n$.

\begin{def.}
For an abelian group $A$, define the chain and cochainquandle complexes

\begin{align}\notag
C_{*}^{Q}(X;A)=C_{*}^{Q}(X)\otimes A, &\quad & \Delta = \partial \otimes id;\\\notag
C_{Q}^{*}(X;A)=Hom(C_{*}^{Q}(X),A), & \quad & \delta = Hom(\partial ,A)
\end{align}

The $n$th quandle homology group and the $n$th quandle cohomology group of a quandle $X$ with coefficient group $A$ are

\begin{align}\notag
H_{n}^{Q}(X;A)=H_n(C_{*}^{Q}(X;A)), &\quad & H_{Q}^{n}(X;A)=H^n(C_{Q}^{*}(X;A))
\end{align}

 We will omit the coefficient group $A$ if $A=\mathbb{Z}$ as usual.
\end{def.}

\begin{rem}
A quandle 2-cocycle $\phi$ satisfies, for a 3-chain $(p,q,r)$:
\begin{align}\notag
0&=(\delta_2(\phi ))(p,q,r)=\phi (\partial_3(p,q,r))=\\\notag
&=\phi (p,r)-\phi (p\triangleright q,r)-\phi (p,q)+\phi (p\triangleright r,q\triangleright r)
\end{align}
i.e.:
\[
\phi (p,r)+\phi (p\triangleright r,q\triangleright r)=\phi (p,q)+\phi (p\triangleright q,r)
\]

Analogously, a 3-cocycle $\theta$ satisfies the relation
\begin{multline}
\theta (p,q,r)+\theta (p\triangleright r,q\triangleright r,s)+\theta (p,r,s)\\
=\theta (p\triangleright q,r,s)+\theta (p,q,s)+\theta (p\triangleright s,q\triangleright s,r\triangleright s)\notag
\end{multline}

It is easy to see by inspection of figures 1,2 and 3 that these express invariance under the type III Reidemeister move and the tetrahedral move, respectively.

Also note that, for any n-cocycle $\phi$
\[
\phi (x_1,...,x_n)=0
\] 
for any n-chain with $x_i=x_{i+1}$, for some $i$. This will express invariance under the type I Reidemeiter move and the type I Roseman moves, for $n=2$ and $n=3$, respectively.
\end{rem}

\begin{figure}
\begin{center}
\mbox{
\epsfxsize=5in
\epsfbox{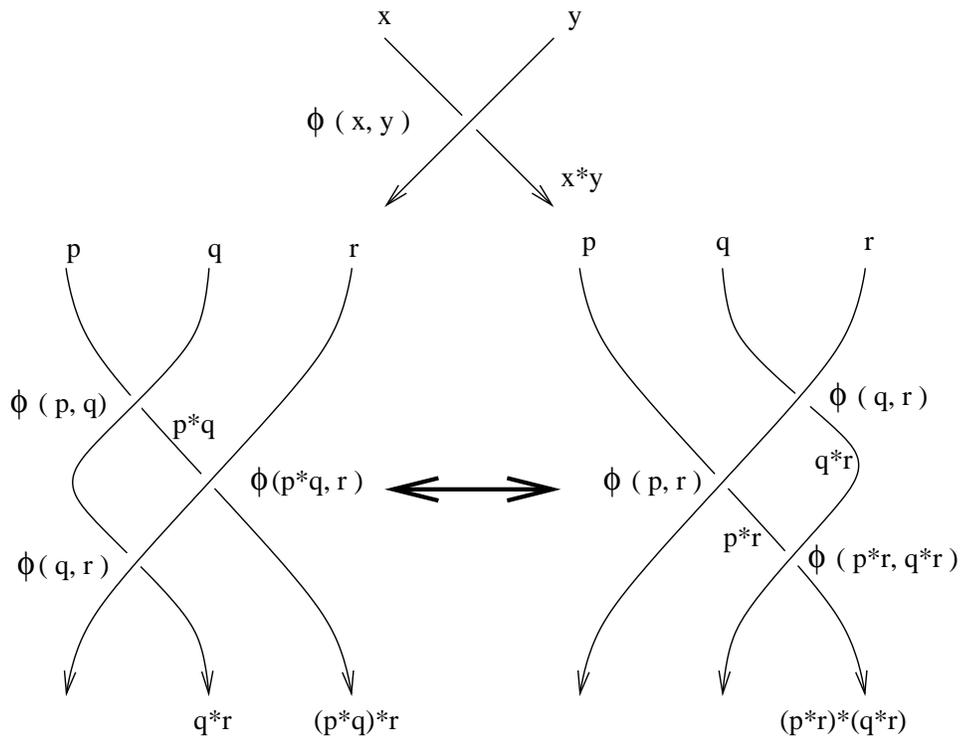}
}
\end{center}
\caption{The 2-cocycle condition and the Reidemeister type III move}
\label{2cocy}
\end{figure}

\begin{figure}
\begin{center}
\mbox{
\epsfxsize=5in
\epsfbox{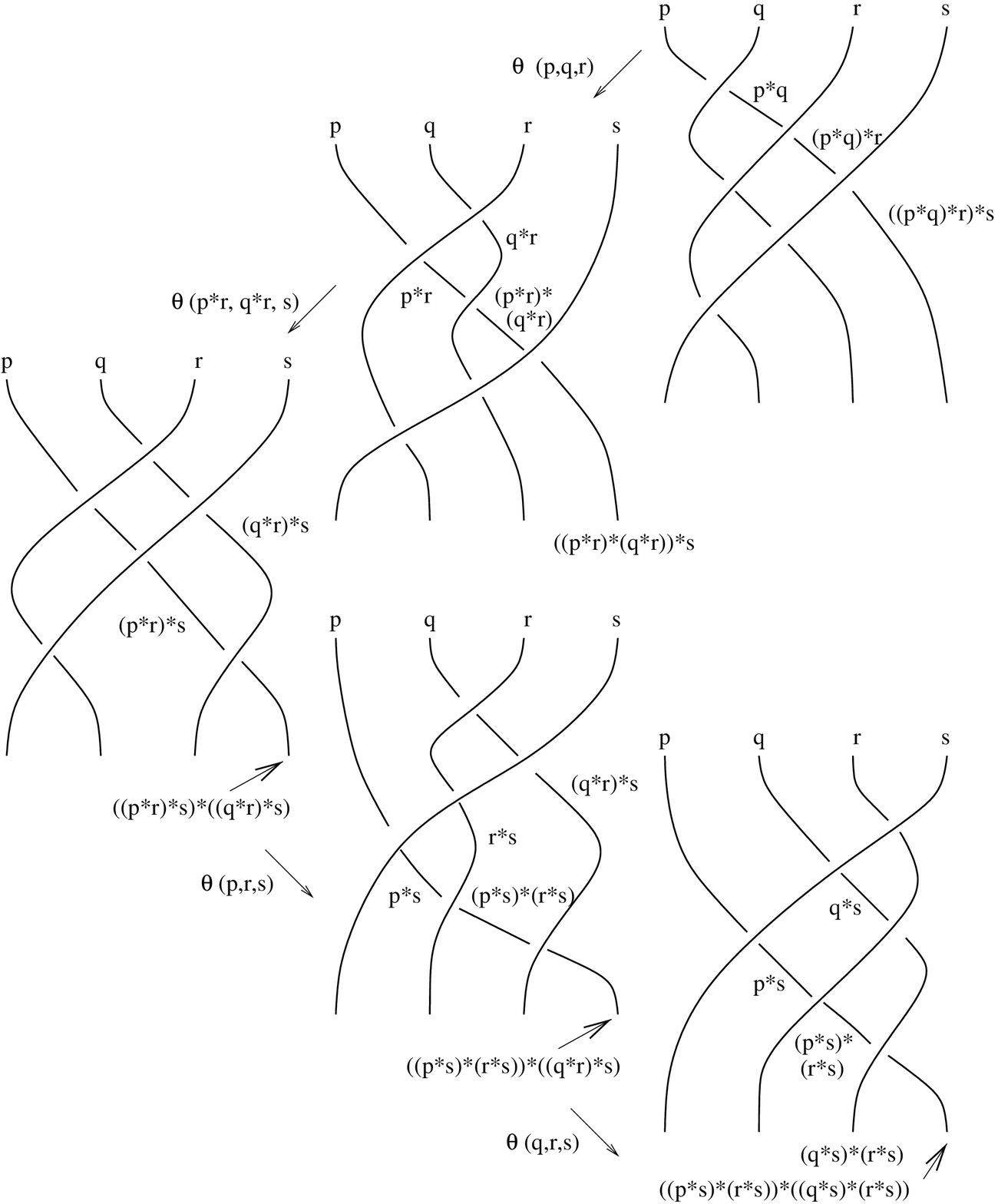}
}
\end{center}
\caption{The tetrahedral move and a cocycle relation, LHS  }
\label{tetraL}
\end{figure}

\begin{figure}
\begin{center}
\mbox{
\epsfxsize=5in
\epsfbox{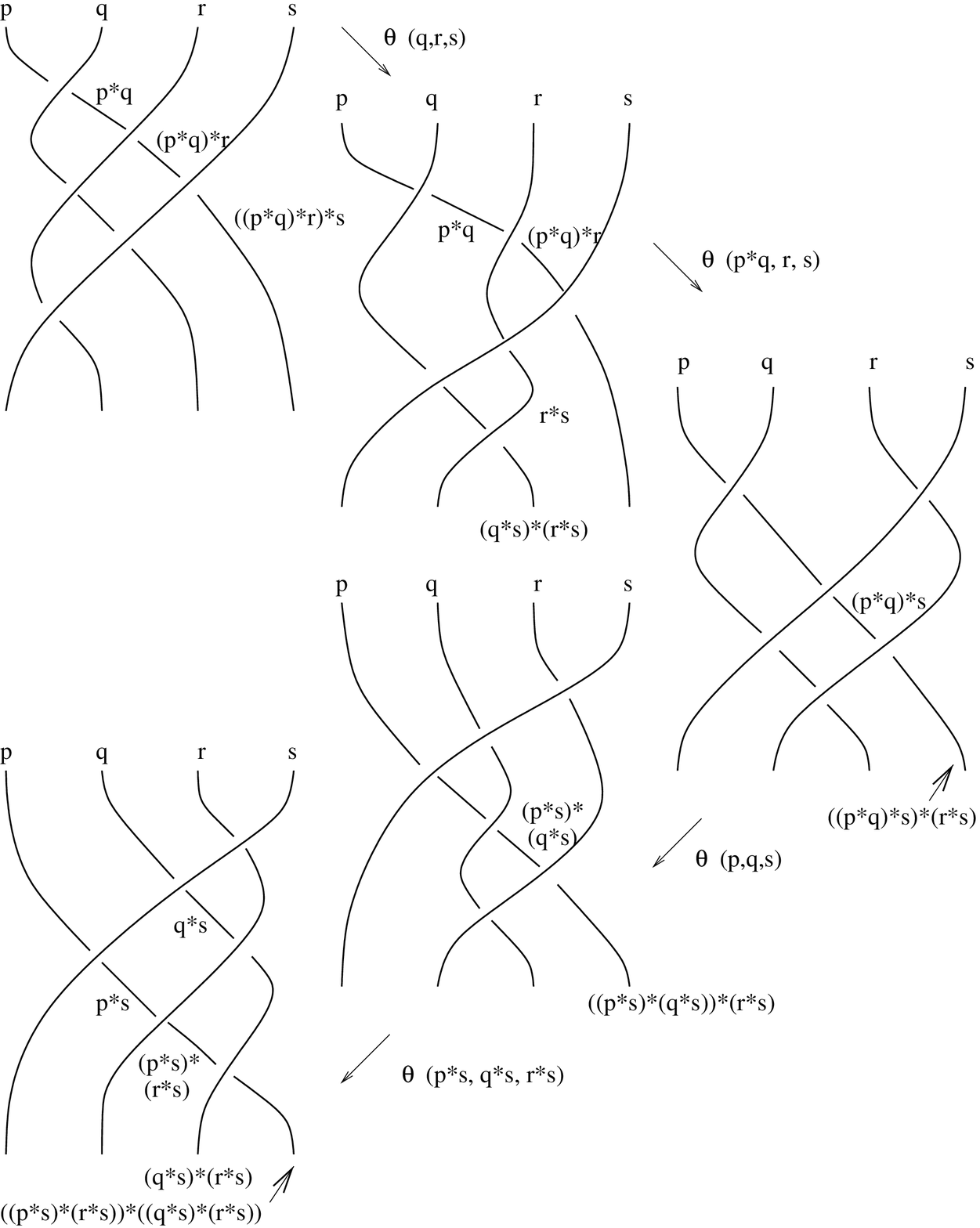}
}
\end{center}
\caption{The tetrahedral move and a cocycle relation, RHS   }
\label{tetraR}
\end{figure}

\section{Classical Knots and Quandle 2-Cocycles} \label{S:Reid}

For the purpose of the discussion in this section we will arbitrarily fix a finite quandle $X$ called the ``labelling quandle'', an abelian group $A$ (denoted multiplicatively), a quandle 2-cocycle $\phi$ (in the quandle cohomology of the labelling quandle with coefficients in $A$), and a knot $K$ (assumed to be oriented); $D$ and $D'$ will stand for diagrams of $K$, identified up to planar isotopy, and $R_D$ and $R_{D'}$ for their sets of arcs.

\begin{def.}[Coloring of a Diagram]A coloring of a diagram, $D$, is a map ${\cal C}:R_D \rightarrow X$ such that, at each crossing:

$$\xy 0;/r1pc/:
,{\xoverv[4.0] |> <{r_1}|{r}>{r_2}}
\endxy$$

if ${\cal C}(r_1)=x$ and ${\cal C}(r)=y$, then ${\cal C}(r_2)=x\triangleright y$. Equivalently, at each crossing:

$$\xy 0;/r1pc/:
,{\xoverv[4.0] |< <{r_1}|{r}>{r_2}}
\endxy$$

if ${\cal C}(r_1)=x$ and ${\cal C}(r)=y$, then ${\cal C}(r_2)=x{\triangleright ^{-1}}y$
\end{def.}
\vskip 5pt

\begin{rem}
Note that colorings are quandle homomorphisms from the fundamental quandle of the knot to the labelling quandle (cf. \cite{dJoyce}).
\end{rem}

\begin{rem}
Note that, locally, the coloring does not depend on the orientation of the under-arc. Also, this notion is well-defined in the sense that, no matter which pair of under-arc and over-arc one colors first, the emergent under-arc will be colored in a consistent way.
\end{rem}

Example:

$$\xy 0;/r1pc/:
,{\xoverv[4.0] |< <{r_1}|{r}>{r_2}}
\endxy$$ 

if ${\cal C}(r)=y$ and ${\cal C}(r_2)=x\triangleright y$, then ${\cal C}(r_1)=(x{\triangleright ^{-1}}y)\triangleright y = x$, which is consistent with the second statement in the previous definition.

\vskip 5pt

\begin{def.}[Set of Colorings of a Knot]We will refer to the set of all colorings of all diagrams of a knot simply as the {\bf set of colorings of a knot}.
\end{def.}
\vskip 5pt

Note that, although there may be lots of different ways of coloring two different diagrams, some of these different colorings are related in a natural way:


\begin{def.}[Relation $\sim$ on the set of colorings of $K$]Let ${\cal C}:R_D \rightarrow X$ and ${\cal C'}:R_{D'} \rightarrow X$ be colorings. We say ${\cal C}\sim {\cal C'}$ if there is a finite sequence of steps taking ${\cal C}$ to ${\cal C'}$ with the details outlined below:
\begin{itemize}

\item If $D'$ is obtained from $D$ by a type I Reidemeister move, then ${\cal C}$ and ${\cal C'}$ are related as follows:
$$\TC\longleftrightarrow\TD$$
with

\[{\cal C}(r)={\cal C'}(r'_1)={\cal C'}(r'_2)
\]
The prescription is analogous for all other cases of type I moves for oriented diagrams. Note that, due to the definition of coloring, the two arcs in the figure on the right have to be assigned the same color, ${\cal C'}(r'_1)={\cal C'}(r'_2)$, for any coloring.

\item If $D'$ is obtained from $D$ by a type II Reidemeister move, then ${\cal C}$ and ${\cal C'}$ are related as follows:
\vskip 10pt
$$\RA\longleftrightarrow\RB$$
\vskip 10pt
\begin{center}
	${\cal C}(r)=x$ \quad and \quad ${\cal C}(s)=y$ \\
	$\sim $\\
	${\cal C'}(r'_1)=x$,\quad ${\cal C'}(s')=y$,\quad ${\cal C'}(r'_2)=x\triangleright y$\quad and\quad ${\cal C'}(r'_3)=(x\triangleright y){\triangleright }^{-1}y=x$
\end{center}
Another case is:
\vskip 10pt
$$\QA\longleftrightarrow\QB$$
\vskip 10pt
\begin{center}
	${\cal C}(r)=x$\quad and\quad ${\cal C}(s)=y$\\
	$\sim $\\
	${\cal C'}(r'_1)=x$,\quad ${\cal C'}(s')=y$,\quad ${\cal C'}(r'_2)=x{\triangleright }^{-1}y$\quad and\quad${\cal C'}(r'_3)=(x{\triangleright }^{-1}y)\triangleright y=x$
\end{center}
Note that the relation between the colorings of the arcs on the right hand side hold for any coloring ${\cal C'}$, from def. 3.1 and thus the passage from ${\cal C'}$ to ${\cal C}$ is well-defined. The remaining cases of Reidemeister II moves between oriented diagrams are dealt with analogously.

\vskip 100pt

 \item For the third Reidemeister move, let us first remark that each line segment can have two possible orientations, giving a total of $2^3=8$ possibilities for each (unoriented) configuration of the three line segments. Since the colorings do not depend on the orientation of the lower line segment then we just have to take care of $2^2=4$ possibilities. Now for each (unoriented) configuration of the lower line segment there are two possibilities for the remaining arcs (one is the top line segment and the other the middle line segment). Rotation through angles of $\frac{2}{3}\pi$ and $\frac{4}{3}\pi$ gives the remaining (unoriented) configurations; but as we have seen before, this is the same as assigning first colors to other arcs. Thus, a priori, there are eight cases to be specified. Further inspection of these eight cases cuts it down to four. We deal here with one of them remarking that the remaining ones are dealt with analogously.

\vskip 10pt

$$\dieR$$
$$\qquad\updownarrow $$
$$\Reid$$

\begin{center}
	${\cal C}(r)=x$,\quad ${\cal C}(s)=y$,\quad and\quad ${\cal C}(t)=z$\quad with\quad ${\cal C}(u)=x{\triangleright }^{-1}z$,\quad ${\cal C}(v)=y{\triangleright }^{-1}z$\quad and\quad ${\cal C}(w)=(x{\triangleright }^{-1}z){\triangleright }^{-1}(y{\triangleright }^{-1}z)=(x{\triangleright }^{-1}y){\triangleright }^{-1}z$\\
	$\sim $\\
	${\cal C'}(r')=x$,\quad ${\cal C'}(s')=y$,\quad and \quad${\cal C'}(t')=z$,\quad ${\cal C'}(u')=x{\triangleright }^{-1}y$,\quad with\quad ${\cal C'}(v')=y{\triangleright }^{-1}z$,\quad and \quad ${\cal C'}(w')=(x{\triangleright }^{-1}y){\triangleright }^{-1}z$
\end{center}

\end{itemize}

\end{def.}

\vskip 20pt

\begin{prop}$\sim$ is an equivalence relation on the set of colorings of $K$.\end{prop}Proof:(i) Obviously, ${\cal C}\sim {\cal C}$;

(ii) Suppose ${\cal C}\sim {\cal C'}$. Reversing the steps that led from ${\cal C}$ to ${\cal C'}$ we get ${\cal C'}\sim {\cal C}$;

(iii) Suppose ${\cal C}\sim {\cal C'}$ ($0\leq n$ moves) and ${\cal C'}\sim {\cal C''}$ ($0\leq m$ moves). Then ${\cal C}\sim {\cal C''}$ ($0\leq n+m$ moves).$\hfill \blacksquare$

\begin{rem}

Suppose ${\cal C}\sim {\cal C'}$ with $D$ and $D'$ related by $n$ Reidemeister moves and realizing the connection between the two colorings as above.

Now suppose there is a different set of $m$ Reidemeister moves that also relate $D$ to $D'$; starting from the coloring  ${\cal C}$ for $D$ what coloring do we get for $D'$ following these latter $m$ moves? Not necessarily  ${\cal C'}$ again. In any case, we get a coloring  ${\cal C''}$ of $D'$ equivalent to the coloring ${\cal C'}$ of $D'$ (cf. proof of proposition above).\end{rem}

\begin{rem}
Let $D$ and $D'$ be two (Reidemeister) equivalent diagrams and fix $n$ Reidemeister moves that take $D$ to $D'$. Now color each of them in all possible ways, i.e. consider the sets $S=\{{\cal C}:R_D \rightarrow X\}$ and $S'=\{{\cal C'}:R_{D'} \rightarrow X\}$. Also, think of the maps:
\[
\begin{array}{clcr}	
	\psi : & S \longrightarrow S'\\
		&{\cal C} \longmapsto {\cal C'}
\end{array}
\quad\txt{ and }\quad
\begin{array}{clcr}
	\psi ': & S' \longrightarrow S\\
		&{\cal C'} \longmapsto {\cal C}
\end{array}
\]
(where, in both cases, ${\cal C} \sim {\cal C'}$ via the above fixed Reidemeister moves). It is readily seen that both maps are surjective and their composite is the identity map. Hence both of them are bijections relating pairs of equivalent colorings and they establish a 1-1 correspondence between $S$ and $S'$.
\end{rem}

\vskip 10pt

In this way, we have just proved that:

\begin{thm}The set of equivalence classes of the colorings associated to a diagram (along with their multiplicities) is a knot invariant.\end{thm}

\vskip 10pt

\begin{cor}The number of equivalence classes of the colorings associated to a diagram is a knot invariant.\end{cor}

\vskip 10pt

\begin{cor}The set of multiplicities referred to in Thm. 3.1 is a knot invariant.\end{cor}

\vskip 10pt

\begin{cor}The total number of colorings associated to a knot diagram is a knot invariant.\end{cor}

\vskip 10pt

Now the problem is that colorings are not easy to deal with. In order to overcome this issue we will introduce the so-called Boltzmann weights which will allow us to come up with yet another but more tractable invariant. 

\vskip 10pt

\begin{def.}[Signs of Crossings]If the pair of the co-orientation of the over-arc and that of the under-arc match the (right-hand) orientation of the plane, then the crossing is called {\bf positive}, otherwise it is {\bf negative}:

\vskip 10pt

\[
\epsilon \biggl(\xy 0;/r1.5pc/:
\xoverv~{(-1,1)}{(1,1)}{(-1,-1)}{(1,-1)}<><>|>>>
\endxy \biggr) =1\qquad 
\epsilon \biggl(\xy 0;/l-1.5pc/:
\xoverv~{(-1,1)}{(1,1)}{(-1,-1)}{(1,-1)}<<|>><\endxy \biggr) =-1\]

\end{def.}

\vskip 20pt

\begin{def.}[$\phi $-assignments at Crossings]Consider a 2-cocycle, $\phi $. Given a coloring ${\cal C}:R_D \rightarrow X$, we assign an element of $A$ to each crossing as follows:

\[
\xy 0;/r1.5pc/:
\xoverv~{(-1,1)}{(1,1)}{(-1,-1)}{(1,-1)}<><>|>>><{r_1}|{r}>{r_2}
\endxy \mapsto \phi (x,y) \qquad \qquad
\xy 0;/l-1.5pc/:
\xoverv~{(-1,1)}{(1,1)}{(-1,-1)}{(1,-1)}<<|>><<{r_1}|{r}>{r_2}\endxy \mapsto \phi (x,y)^{-1}\]
where as before, ${\cal C}(r_1)=x$ and ${\cal C}(r)=y$.

\vskip 20pt

We will also write ${\phi_{\tau}}^{\epsilon_{\tau}}$ where $\tau$ denotes the crossing where the assignment is being made.
\end{def.}

\begin{def.}[Boltzmann weight]Given a 2-cocycle $\phi $ and a coloring ${\cal C}:R_D \rightarrow X$ the Boltzmann weight associated to ${\cal C}$ is:
\[
B({\cal C}) = \prod _{\tau } {\phi _{\tau }}^{\epsilon_{\tau }}
\]

where $\tau $ runs over the set of crossings of $D$.

\end{def.}
\vskip 20pt

Note that our notion of Boltzmann weight differs from the one in \cite{jsCetal} (where Boltzmann weights are the individual factors in the above product) and is closer, in meaning, to the statistical mechanical one (cf. \cite{hbCallen}).

\vskip 20pt

\begin{prop}${\cal C}\sim {\cal C'}$ implies $B({\cal C})=B({\cal C'})$.\end{prop}Proof: It is enough to show the statement for the moves on colorings corresponding to the Reidemeister moves.

(i) For type I moves the addition of a crossing $\tau$ with both arcs colored with the same quandle element $x$ contributes a factor ${\phi_{\tau}}^{\epsilon_{\tau}}=1^{\pm 1}$ to the Boltzmann weight.

(ii) For type II moves one has the addition of two crossings $\tau_1$ and $\tau_2$ with $\epsilon_{\tau _1}+\epsilon_{\tau _2}=0$ which contribute a factor ${\phi _{\tau _1}}^{\epsilon_{\tau _1}} \cdot  {\phi _{\tau _2}}^{\epsilon_{\tau _2}}=1$ to the Boltzmann weight since $\phi _{\tau _1}=\phi _{\tau _2}$.

(iii) The defining equation for $\phi$ implies invariance under a type III Reidemeister move with a particular set of signs in the crossings (cf. Remark 2.1 and \cite{jsCetal}); invariance under type III moves with other sets of signs follows from particular combinations of the former type III move with type II moves (cf. \cite{lhKauffman}).$\hfill \blacksquare$

\vskip 20pt

Hence, the Boltzmann weights are constant on $\sim$-equivalence classes.

The following corollaries are trivial consequences of Prop. 3.3

\vskip 20pt

\begin{cor}The set of Boltzmann weights associated to the colorings of a knot diagram (along with their multiplicities) is a knot invariant.
\end{cor}Proof: (Reidemeister) equivalent diagrams have their colorings paired with equivalent colorings in each pair (cf. Remark 3.3). Since equivalent colorings are assigned the same Boltzmann weight, the same set of Boltzmann weights is assigned to each of the diagrams.$\hfill \blacksquare$

\vskip 20pt

\begin{cor}The set of multiplicities referred to in cor. 3.4 is a knot invariant.
\end{cor}

\vskip 20pt

In this way we recover the result of \cite{jsCetal}:

\begin{cor}$\sum_{\cal C}B(\cal C)$ is a knot invariant (where $\cal C$ runs over the set of all colorings).
\end{cor}

\vskip 20pt

\begin{cor}Any symmetric function on all Boltzmann weights is a knot invariant.
\end{cor}

\vskip 20pt

\begin{prop}If $\phi$ is a coboundary, $\sum_{\cal C}B(\cal C)$ is the number of colorings.
\end{prop}Proof: Cf. \cite{jsCetal}.$\hfill \blacksquare$

\vskip 20pt

\begin{rem}
Note that proposition 3.4 and corollary 3.6 imply that the number of colorings is a knot invariant, which was already shown (cf. cor. 3.3).
\end{rem}

\vskip 20pt

\begin{prop}Cohomologous cocycles give rise to the same set of Boltzmann weights.
\end{prop}Proof: Cf. \cite{jsCetal}.$\hfill \blacksquare$

\vskip 20pt

Now, suppose $A$ is a (multiplicative) subgroup of $\mathbb{R}^{+} -\{0\}$.

\begin{prop}For any real number, $ \nu$, $\sum_{\cal C}B(\cal C)^{\nu}$ is a knot invariant.
\end{prop}Proof: Trivial.$\hfill \blacksquare$

\vskip 20pt

Further, there exists a map $\psi$ such that, for any diagram $D$:

\[
\phi_{\tau}=\exp \psi_{\tau}
\]
 for any crossing, $\tau $ of $D$. Hence,
\[
B({\cal C }) = \prod _{\tau } {\phi _{\tau }}^{\epsilon_{\tau }} = \prod _{\tau } {\exp (\psi _{\tau })}^{\epsilon_{\tau }} =  \prod _{\tau } e^{{\epsilon_{\tau }}\psi _{\tau }} =  e^{\sum _{\tau }{\epsilon_{\tau }}\psi _{\tau }}
\]
So, in the notation of the previous proposition, taking $\nu =-\frac{1}{kT}$, we have the "partition function":

\[
Z(T) = \sum _{\cal C} B({\cal C})^{-\frac{1}{kT}}=\sum _{\cal C} {e^{{-\frac{1}{kT} }{\sum _{\tau }\epsilon_{\tau }\psi _{\tau }}}}
\]

So, in this way we introduced the temperature parameter. Because the number of colorings is constant one would say that there is, also, a volume parameter present, here. This establishes an interesting parallel with the canonical ensemble in statistical mechanics (cf. \cite{hbCallen}).

\vskip 20pt

We now have a new of calculating invariants of a knot using quandle cohomology. We obtain a 1-parameter family of invariants $Z(T)$, instead of a single number.

\section{Knotted Surfaces and Quandle 3-Cocycles} \label{S:Roseman}

We briefly recall the notion of a knotted surface diagram (cf. \cite{CetS}). Let $f : F \mapsto \mathbb{R}^4$ denote a smooth embedding of a closed surface $F$ into 4-dimensional space. By deforming the map $f$ slightly by an ambient isotopy of $\mathbb{R}^4$, if necessary, we may assume that $p\circ f$ is a general position map, where $p : \mathbb{R}^4 \mapsto \mathbb{R}^3$ denotes the orthogonal projection onto an affine subspace (cf. \cite{dRoseman2}).

Along the double curves, one of the sheets (called the over-sheet) lies farther than the other (the under-sheet) with respect to the projection direction. The under-sheets are consistently broken in the projection, and such broken surfaces are called knotted surface diagrams (cf. \cite{CetS}).

When the surface is oriented, we take normal vectors $\vec{n}$ to the projection of the surface such that the triple $(\vec{v}_1,\vec{v}_2,\vec{n})$ matches the orientation of 3-space, where $(\vec{v}_1,\vec{v}_2)$ defines the orientation of the surface. Such normal vectors are defined on the projection at all points other than the branch points.

So, for the purpose of the discussion in this section we will arbitrarily fix a finite quandle $X$ (the labelling quandle), an abelian group $A$ (denoted multiplicatively), a quandle 3-cocycle $\phi$ (in the quandle cohomology of the labelling quandle with coefficients in $A$), an orientable knotted surface $K$, together with two of its (broken surface) diagrams, $D$ and $D'$ (identified up to isotopy in $\mathbb{R}^3$ and such that the composite of the embedding with the projection is in general position) and where $R_D$ and $R_{D'}$ stand for the regions of the projection of the surface bounded by the double lines, i.e., the regions one is left with when the projection of the surface is consistently broken along the double lines in order to become a broken surface diagram (cf. \cite{CetS}). We will, further, assume that the surface is oriented and hence, each of the above regions is assigned a normal as discussed in the previous paragraph.

We now present the analogs of the definitions and propositions of the previous section.

\begin{def.}[Coloring of a diagram]A coloring of a knotted surface diagram, $D$, is a map, ${\cal C}:R_D \rightarrow X$, satisfying the following condition.

At a double point curve two coordinate planes intersect locally. One is the over-sheet $r$, the other is the under-sheet which is broken into two components, say $r_1$ and $r_2$. A normal of the over-sheet points to one of the components, say $r_2$. If ${\cal C}(r_1)=x$ and ${\cal C}(r)=y$, then we require that ${\cal C}(r_2)=x\triangleright y$; if the normal of the over-sheet does {\bf not} point to $r_2$, then ${\cal C}(r_2)=x\triangleright ^{-1}y$.
\end{def.}
\vskip 20pt

\begin{prop}The above condition is compatible at each triple point.\end{prop}Proof: See \cite{jsCetal}.$\hfill \blacksquare$

\vskip 20pt

The treatment of knotted surfaces was boosted by the contribution of  D. Roseman who proved that diagrams in 3-space of equivalent embeddings of surfaces in 4-space are related by moves (in pretty much the same way that planar diagrams of (classical) knots are related by the so-called Reidemeister moves). These moves are now known as Roseman moves (cf. \cite{dRoseman1} and also \cite{CetS}).


\begin{figure}
\begin{center}
\mbox{
\epsfxsize=5in
\epsfbox{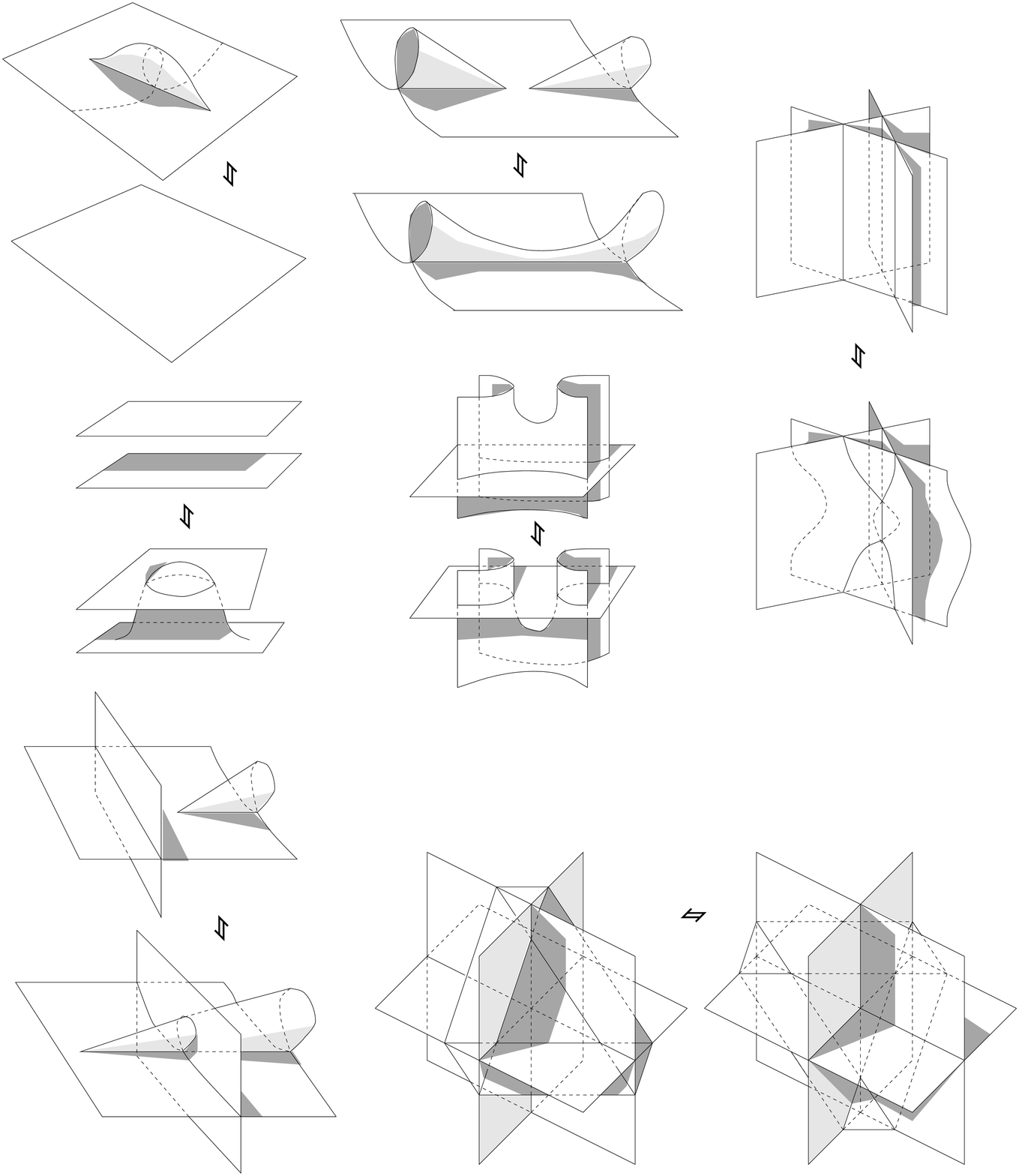}
}
\end{center}
\caption{Roseman moves for knotted surface diagrams  }
\label{rose}
\end{figure}


Following \cite{CetS}, page 43, we name these different Roseman moves as follows (see figure 4):

The first two moves depicted on the top row of the figure are both called type-I/type-I-inverse moves; the move in the left column is called a type-I bubble move; the move in the second column is called a type-I saddle move. The two moves on the second row are called type-II/type-II-inverse moves; the move on the left is called a type-II bubble move; the move on the right is called a type-II saddle move. The remaining are type-III/type-III-inverse moves. In the move in the far right column on the top a pair of oppositely signed triple points is either cancelled or introduced. The move depicted on the left of the bottom row is called ``passing a branch point through a sheet''. The move on the bottom right is called the tetrahedral move or the quadruple point move.

Note that those should have been broken surface diagrams instead but the excess of detail would blur the interpretation of the moves. Thus each diagram represents a number of different cases corresponding to different choices of over and under sheets. Moreover the over/under assignment is to some extent irrelevant for type I and type II moves since changing it corresponds to either interchanging the sheets, or reversing the passage, or both.

\vskip 20pt

\begin{def.}[Set of Colorings of a Knotted Surface]We will refer to the set of all colorings of all diagrams of a knotted surface simply as the {\bf set of colorings of a knotted surface}.
\end{def.}
\vskip 20pt

\begin{def.}[Relation $\sim$ on the set of colorings of $K$]Let ${\cal C}:R_D \rightarrow X$ and ${\cal C'}:R_{D'} \rightarrow X$ be colorings. Say ${\cal C}\sim {\cal C'}$ if there is a finite sequence of steps taking ${\cal C}$ to ${\cal C'}$ with the details outlined below (the descriptions are based on figure 4):
\begin{itemize}

\item Suppose $D$ and $D'$ are related by a type I Roseman move.

For the bubble move, since there is only one  region at issue on either side of the move, say $r$ and $r'$, that region is colored with the same quandle element for both sides of the move:

\[
{\cal C}(r)={\cal C'}(r')
\]

The saddle move is the passage from a local configuration of one region, say $r$, to a local configuration of two regions, say $r'_1$ and $r'_2$, with a double-line where one of the regions is both the over-sheet and the under-sheet with respect to that double-line. In this way:
\begin{center}
	${\cal C}(r_1)=x$\\
	$\sim $\\
	${\cal C'}(r'_1)={\cal C'}(r'_2)=x$
\end{center}
Note that for the two regions configuration, since one of the regions is both the under-sheet and the over-sheet, quandle axioms I and II imply that both regions have to have the same color.

\item Suppose $D$ and $D'$ are related by a type II Roseman move, i.e. the II bubble or II saddle. We will assume that whenever a horizontal plane intersects the rest of the surface it is the top sheet. Making the opposite choice corresponds to a deformation plus rotation of $\pi$ (bubble move) or a deformation, a rotation of $\frac{\pi}{2}$ around an axis perpendicular to the horizontal plane plus exchanging the two sides of the move (saddle move).

Calling the two sheets of the upper diagram of the type II bubble move $r_1$ and $r_2$ (with $r_1$ the lower sheet in the diagram) and calling the horizontal over-sheet of the lower diagram $r'$, the lower part of the curved under-sheet $r'_1$ and the upper part of the curved under-sheet $r'_2$ we have:
\begin{center}
	${\cal C}(r_1)=x$\quad and\quad ${\cal C}(r_2)=y$\\
	$\sim$\\
	${\cal C'}(r'_1)=x$,\quad ${\cal C'}(r')=y$\quad  and\quad ${\cal C'}(r'_2)=x\triangleright^{\pm 1}y$
\end{center}
where the coloring of $r'_2$ depends on the orientation of the horizontal sheet.

Now the saddle move. For the upper diagram, let us call $r$ the horizontal sheet, $r_1$ the front lower part of the under-sheet, $r_3$ the rear lower part of the under-sheet and $r_2$ the upper part of the under-sheet. For the lower diagram, let us call $r'$ the over-sheet, $r'_1$ the left upper part of the under-sheet, $r'_3$ the right upper part of the under-sheet and $r'_2$ the lower part of the under-sheet. In this way:
\begin{center}
	${\cal C}(r)=y$\quad and\quad ${\cal C}(r_1)=x={\cal C}(r_3)$\\
	$\sim$\\
	${\cal C}(r_2)=x\triangleright^{\pm 1}y$,\quad ${\cal C'}(r)=y$, ${\cal C'}(r'_2)=x$\quad and \quad${\cal C'}(r'_1)={\cal C'}(r'_3)=x\triangleright^{\pm 1}y$
\end{center}
The $\pm 1$ reflect the dependence on the orientation.

\vskip 15pt

 \item Suppose $D$ and $D'$ are related by a type III Roseman move. 

Consider the "branch point through a sheet" move. This involves two planes say, the horizontal and the vertical planes. The horizontal plane has a branch point attached to it which is a folding of the plane about a point. In the upper diagram the vertical plane intersects the horizontal plane along a line not involving either the branch point or the folding that produces it; in the lower diagram the vertical plane intersects the folding  of the horizontal plane (but not the branch point). Let us first assume that the horizontal plane is the top sheet in both diagrams. For the upper diagram we will call the horizontal plane $r$, the lower region of the vertical plane $r_1$ and the upper region $r_2$, and, for the lower diagram, we will call the horizontal plane $r'$, the lower part of the vertical plane $r'_1$, the upper part of the vertical plane $r'_2$ and the disk on the vertical plane bounded by the folding $r'_3$. We have:
\begin{center}
	${\cal C}(r_1)=x$,\quad ${\cal C}(r)=y$\quad and\quad ${\cal C}(r_2)=x\triangleright^{\pm 1}y$\\
	$\sim$\\
	${\cal C'}(r'_1)=x$,\quad ${\cal C'}(r')=y$,\quad ${\cal C'}(r'_2)=x\triangleright^{\pm 1}y$\quad and \quad ${\cal C'}(r'_3)=x$
\end{center}
Now, assume the horizontal plane is bottom sheet. In the upper diagram let us call the vertical plane $r$, the left part of the horizontal plane $r_1$ and the right part $r_2$; in the lower diagram, let us call the vertical plane $r'$, the left part of the horizontal plane and the right part of the horizontal plane $r'_1$ and $r'_2$, respectively-. In this way:
\begin{center}
	${\cal C}(r_1)=x$,\quad ${\cal C}(r)=y$\quad and\quad ${\cal C}(r_2)=x\triangleright^{\pm 1}y$\\
	$\sim$\\
	${\cal C'}(r'_1)=x$,\quad ${\cal C'}(r')=y$,\quad and\quad ${\cal C'}(r'_2)=x\triangleright^{\pm 1}y$
\end{center}
The $\pm 1$ reflects the dependence on orientation.

Consider the creation/cancellation of a pair of triple points move. There are three cases to be considered: either the plane that folds about the others is top sheet or middle sheet or bottom sheet. The remaining planes play a symmetric role so there is no point in considering interchange of their positions. The plane that subsequently folds about the others will be referred to as the folding sheet; the other planes will be referred to as the right plane (the one which is more to the right in the figure) and the left plane (analogously). We will take the right plane always to be below the left plane in the embedding. So let us first assume that the folding sheet is the top sheet. In the upper diagram let us call the folding sheet $r$, the leftmost part and the rightmost part of the left plane $s_1$ and $s_2$, respectively, and the leftmost part, the middle part and the rightmost part of the right plane $t_1$, $t_2$, $t_3$, respectively. In the lower diagram the folding sheet, $r'_1$, splits the left plane again into two parts, $s'_1$ (the leftmost) and $s'_2$ (the rightmost), and the right plane into five parts. The previous leftmost part of it is split into a new leftmost part $t'_1$ and a righmost part $t'_2$; the remainder of the right plane is split into a top part $t'_3$, a new rightmost part $t'_4$ and a bottom part $t'_5$. In this way:
\begin{center}
	${\cal C}(s_1)=x$,\quad ${\cal C}(t_1)=y$,\quad ${\cal C}(r)=z$ ,\quad ${\cal C}(s_2)=x\triangleright^{\epsilon_1}z$,\quad ${\cal C}(t_2)=y\triangleright^{\epsilon_2}(x\triangleright^{\epsilon_1}z)$,\quad ${\cal C}(t_3)=[y\triangleright^{\epsilon_2}(x\triangleright^{\epsilon_1}z)]\triangleright^{-\epsilon_1}z=(y\triangleright^{-\epsilon_1}z)\triangleright^{\epsilon_2}x$\\
	$\sim$\\
	${\cal C'}(s'_1)=x$,\quad ${\cal C'}(t'_1)=y$,\quad ${\cal C'}(r')=z$,\quad ${\cal C'}(s'_2)=x\triangleright^{\epsilon_1}z$,\quad ${\cal C'}(t'_2)=y\triangleright^{-\epsilon_1}z$,\quad ${\cal C'}(t'_3)={\cal C'}(t'_5)=y\triangleright^{\epsilon_2}(x\triangleright^{\epsilon_1}z)$,\quad ${\cal C'}(t'_4)=(y\triangleright^{-\epsilon_1}z)\triangleright^{\epsilon_2}x$
\end{center}
The $\epsilon$'s reflect the dependence on the orientation.

Let us now assume that the folding sheet is the middle sheet. In the upper diagram let us call the leftmost and the rightmost parts of the folding plane $s_1$ and $s_2$, respectively, the left plane $r$ and the leftmost part, the middle part and the rightmost part of the right plane $t_1$, $t_2$, $t_3$, respectively. In the lower diagram, the left plane, $r'$, together with the leftmost and rightmost parts of the folding sheet, $s'_1$ and $s'_2$, respectively, split the right plane into five parts. As before the previous leftmost part of it is split into a new leftmost part $t'_1$ and a righmost part $t'_2$; the remainder of the right plane is split into a top part $t'_3$, a new rightmost part $t'_4$ and a bottom part, $t'_5$.We have:
\begin{center}
	${\cal C}(s_1)=x$,\quad ${\cal C}(t_1)=y$,\quad  ${\cal C}(r)=z$,\quad  ${\cal C}(s_2)=x\triangleright^{\epsilon_1}z$,\quad  ${\cal C}(t_2)=y\triangleright^{\epsilon_1}z$,\quad and \quad ${\cal C}(t_3)=(y\triangleright^{\epsilon_1}z)\triangleright^{\epsilon_2}(x\triangleright^{\epsilon_1}z)=(y\triangleright^{\epsilon_2}x)\triangleright^{\epsilon_1}z$\\
	$\sim$\\
	${\cal C'}(s'_1)=x$,\quad ${\cal C'}(t'_1)=y$,\quad ${\cal C'}(r')=z$,\quad ${\cal C'}(s'_2)=x\triangleright^{\epsilon_1}z$,\quad ${\cal C'}(t'_2)=y\triangleright^{\epsilon_2}x$,\quad ${\cal C'}(t'_3)={\cal C'}(t'_5)=y\triangleright^{\epsilon_1}z$,\quad and \quad ${\cal C'}(t'_4)=(y\triangleright^{\epsilon_2}x)\triangleright^{\epsilon_1}z$\\
\end{center}
Again the $\epsilon$'s reflect the dependence on the orientation.

Finally, let us assume that the folding sheet is bottom sheet. In the upper diagram we call the left plane $r$, the leftmost and rightmost parts of the right plane $s_1$ and $s_2$, respectively, and the leftmost, middle and rightmost parts of the folding sheet $t_1$, $t_2$ and $t_3$, respectively. In the lower diagram the folding sheet is split into a leftmost part, a rightmost part and three central parts, top, middle and bottom, respectively. We call the left plane $r'$, the leftmost and rightmost parts of the right plane $s'_1$ and $s'_2$, respectively, and the leftmost and rightmost parts of the folding sheet $t'_1$ and $t'_2$, respectively, and the top, middle and bottom parts of the folding sheet, $t'_3$, $t'_4$ and $t'_5$, respectively. In this way:
\begin{center}
	${\cal C}(t_1)=x$,\quad ${\cal C}(s_1)=y$,\quad ${\cal C}(r)=z$,\quad ${\cal C}(s_2)=y\triangleright^{\epsilon_1}z$,\quad ${\cal C}(t_2)=x\triangleright^{\epsilon_1}z$,\quad and \quad ${\cal C}(t_3)=(x\triangleright^{\epsilon_1}z)\triangleright^{\epsilon_2}(y\triangleright^{\epsilon_1}z)=(x\triangleright^{\epsilon_2}y)\triangleright^{\epsilon_1}z$\\
	$\sim$\\
	${\cal C'}(t'_1)=x$,\quad ${\cal C'}(s'_1)=y$,\quad ${\cal C'}(r')=z$,\quad ${\cal C'}(s'_2)=y\triangleright^{\epsilon_1}z$,\quad ${\cal C'}(t'_3)={\cal C'}(t'_5)=x\triangleright^{\epsilon_1}z$, ${\cal C'}(t'_4)=x\triangleright^{\epsilon_2}y$\quad and \quad ${\cal C'}(t'_2)=(x\triangleright^{\epsilon_2}y)\triangleright^{\epsilon_1}z$
\end{center}
As before the $\epsilon$'s denote dependence on orientation.

Finally, consider the tetrahedral move. The specifications of the colorings are given in figures 2 and 3 for a particular case ( note that the technique for representing embeddings here is the so-called movie moves - cf. \cite{CetS}). The remaining cases are dealt with analogously.

\end{itemize}
\end{def.}

\vskip 20pt

\begin{prop}
$\sim$ is an equivalence relation on the set of colorings of $K$.\end{prop}Proof: Analogous to Prop. 3.1 $\hfill \blacksquare$

\vskip 10pt

The analogs of Remarks 3.3 and 3.4 are also valid here. Hence, we have:

\begin{thm}The set of equivalence classes of the colorings associated to a diagram (along with their multiplicities) is an invariant for knotted surfaces.\end{thm}

\vskip 10pt

The analogous corollaries to cors. 3.1-3.3 also hold here.
\vskip 10pt

As we did before, we will now introduce the Boltzmann weights.

\vskip 10pt

\begin{def.}[Signs of Triple Points]Note that when three sheets form a triple point, they have relative positions top, middle, bottom with respect to the projection direction of $p : \mathbb{R}^4\mapsto \mathbb{R}^3$. The sign of a triple point is positive if the normals of top, middle, bottom sheets in this order match the orientation of the 3-space. Otherwise the sign is negative. We denote the sign of a triple point, $\tau$, by $\epsilon_{\tau}$. We use the right-hand rule convention for the orientation of 3-space.
\end{def.}
\vskip 20pt

\begin{def.}[$\phi $-assignments at Triple Points]Let $\phi $ be a three cocycle and recall that it satisfies

\begin{multline}\label{E:mm3}
\phi (p,q,r)+\phi (p\triangleright r,q\triangleright r,s)+\phi (p,r,s)\\
=\phi (p\triangleright q,r,s)+\phi (p,q,s)+\phi (p\triangleright s,q\triangleright s,r\triangleright s)\notag
\end{multline}
At each triple point, $\tau$, the $\phi $ function will be assigned as follows:
 Let O be the octant from which all normal vectors of the three sheets point outwards; let $x$, $y$, $z$ be the colors of the bottom, middle and top sheets, respectively, that bound the region O. Let $\epsilon_{\tau}$ be the sign of the triple point. Then the assignment at this triple point is $\phi (x,y,z)^{\epsilon_{\tau}}$.
\end{def.}
\vskip 40pt

\begin{def.}[Boltzmann weight]Given a coloring ${\cal C}:R_D \rightarrow X$ the Boltzmann weight associated to it is:
\[
B({\cal C}) = \prod _{\tau } {\phi _{\tau }}^{\epsilon_{\tau }}
\]

where $\tau $ runs over the set of triple points of $D$.
\end{def.}

\vskip 20pt

\begin{prop}
${\cal C}\sim {\cal C'}$ implies $B({\cal C})=B({\cal C'})$.\end{prop}Proof: We are here mimicking the proof of Thm. 5.6 in \cite{jsCetal}. The Boltzmann weights depend only on triple points. Thus, any move not involving triple points will not affect the Boltzmann weights. In this way we just have to check what happens when the underlying diagrams are related by Roseman moves with triple points. These are: (1) the creation or cancellation of a pair of oppositely signed triple points; (2) moving a branch point through a sheet; (3) the tetrahedral move.

In the first case, the pair of triple points have opposite signs, so for a given coloring, the two contributing factors of the state-sum cancel. In the second case, the branch point occurs on either the bottom/middle sheet or on the top/middle sheet, and these sheets have the same color. Since the weighting of the proximate triple point is a quandle cocycle (so $\phi(x,x,y)=\phi(x,y,y)=1$), this factor does not contribute to the state-sum.

In the third case, there are several possible tetrahedral moves to consider that depend on (a) the local orientation of the sheets around the tetrahedron, and (b) the signs of the triple points that are the vertices of the tetrahedron. The definition of the 3-cocycle implies that the state-sum is invariant under one of these possible choices. We will move a given tetrahedral move so that the planes involved coincide with planes in this standard position, but have possibly differing crossings or orientations. Then a generalization of Turaev's technique (cf. \cite{vTur1} and \cite{vTur2}) to dimension 4 shows that the given move follows from the fixed move and invariance under adding or subtracting a cancelling pair of triple points (see \cite{jsCetal}) which finishes the proof.$\hfill \blacksquare$

\vskip 20pt

The results of the previous section now go through in an entirely analogous fashion. In particular we recover the result in \cite{jsCetal} that $\sum_{\cal C}B(\cal C)$ is an invariant for knotted surfaces as well as the stronger result that the set of Boltzmann weights, the set of their multiplicities and any symmetric function of the Boltzmann weights are invariants. Again one such invariant is the finite temperature partition function

\[
Z(T) = \sum _{\cal C} B({\cal C})^{-\frac{1}{kT}}=\sum _{\cal C} {e^{{-\frac{1}{kT} }{\sum _{\tau }\epsilon_{\tau }\psi _{\tau }}}}
\]

where $\phi_{\tau}=\exp \psi_{\tau}$ and $T$ is the temperature, suggesting a statistical mechanical interpretation.

\section{Future Directions} \label{S:concl}

In future work we will study particular examples for specific quandles and address problems which could be solved by these methods as well as explore the connections with Statistical Mechanics.

\end{document}